\documentclass[smallextended]{svjour3}       % onecolumn (second format)
\smartqed  % flush right qed marks, e.g. at end of proof
\usepackage{amsfonts}
\usepackage{amsmath}
\usepackage{graphicx}
\usepackage{epsfig}
\usepackage{amssymb}
\usepackage{tikz}
\usetikzlibrary{positioning}
\usetikzlibrary{shapes}

\usepackage{comment}

\newenvironment{PrfFact}{{\bf Proof.}}{{\hfill{$\blacksquare$\\ }}}

\newlength\flitemwidth

\journalname{Journal of Combinatorial Optimization}

\begin{document}

\title{The Stein theorem for loopless $2$-connected plane multigraphs}

\author{Jan~Florek}

\institute{J.~Florek \at  Faculty of Pure and Applied Mathematics,
 Wroclaw University of Science and Technology,
 Wybrze\.{z}e Wyspia\'nskiego 27,
50--370 Wroc{\l}aw, Poland\\
\email{jan.florek@pwr.edu.pl}}

\date{Received: date / Accepted: date}

\maketitle

\begin{abstract}
Stein proved that for each  simple plane triangulation $H$  there exists  a  partitioning of the vertex set of $H$ into two subsets each of which induces a forest if and only if the dual $H^{*}$ has a Hamilton cycle.We extend the Stein theorem for graphs in the family of all loopless $2$-connected plane multigraphs and we prove some other equivalent results.
%\begin{theorem}\label{theorem1} 
%Let $G \in {\cal G}$ and suppose that $S$ and $T$ are disjoint induced subgraphs of $G$ which together contain the vertex set of $G$. The following conditions are equivalent:

%1. $S$ and $T$ are acyclic and $S\cap C$  is a path or a vertex  for every facial cycle~$C$ of $G$,

%2. $ \{e^{*} \in E(G^{*}) : e \in  E[S, T] \}$ is  the edge set of a Hamilton cycle in $G^{*}$, 

%3. $S$ and $T$ are trees,

%4. $S$ is a  tree and $S\cap C$ is a path or a vertex for every facial cycle $C$ of $G$.
%\end{theorem}
%\footnotetext{2010 \textit{Mathematics Subject Classification}: 05C45, 05C10.}
%\footnotetext{\textit{Key words and phrases}: Barnette's conjecture, Hamilton cycle, induced tree, contraction of a path.}

 \end{abstract}


%\end{frontmatter}
\section{Introduction}
We use \cite{flobar1} as reference for undefined terms.
 
Stein  (\cite{flobar4}) expresses hamiltonicity in terms of the dual graph. Let $H$ be a simple plane triangulation. Stein proved that there exists  a  partitioning of the vertex set of $H$ into two subsets each of which induces a forest if and only if the dual $H^{*}$ has a Hamilton cycle. Let $K$ be a connected plane graph. Hakimi and Schmeichel \cite{flobar2}) proved that there exists  a  partitioning of the vertex set of~$K$ into two subsets each of which induces a forest if and only if the dual $K^{*}$ has an eulerian spanning subgraph. Skupie\'{n} \cite{flobar3} proved the Hakimi and Schmeichel theorem for  connected plane multigraphs without loops.

Let $\cal G$ be the family of all loopless $2$-connected plane multigraph. Each face of a graph belonging to ${\cal G}$ is bounded by a cycle, called a \textsl{facial cycle} (Bondy and  Murty \cite{flobar1} Theorem 10.7). Suppose that $G \in \cal G$, $V(G)$  is the vertex set and $E(G)$ is the edge set of $G$.  If $\{X, Y\}$  is a partition of $V(G)$, the set $E[X, Y] = \{xy \in E(G): x \in X \hbox{ and }  y\in Y\}$ is called an \textsl{edge cut} of $G$. A \textsl{bond} of $G$ is a minimal nonempty edge cut of $G$. Each edge cut contains a bond (Bondy and  Murty \cite{flobar1} Theorem 2.14). If $E[X, Y]$ is a bond, then both $G[X]$ and $G[Y]$ are connected, because $G$ is connected (Bondy and  Murty \cite{flobar1} Theorem 2.15). Let $G^{*}$ denote the dual graph of $G$. If $S$ is a proper subgraph of $G$, then we denote $S^{*} =  \{e^{*}\in E(G^{*}) :  e \in E(S) \}$. There exists a relationship between cycles of $G$ and bonds of  $G^{*}$.  Namely, $C$ is a cycle in $G$ if and only if $C^{*}$  is a bond of $G^{*}$ (Bondy and  Murty \cite{flobar1} Theorem 10.16). %If $C$ is a Hamilton cycle in $G$, then both components of $G^{*} \backslash C^{*}$ are trees (see Bondy and  Murty \cite{flobar1} 10.2.11). 
 Notice that the dual of a graph in ${\cal G}$ belongs to ${\cal G}$ (see Bondy and  Murty \cite{flobar1} Theorem 10.14).  

We prove the following theorem which extend the Stein theorem for graphs in~$\cal G$.
\begin{theorem}\label{theorem1} 
Let $G \in {\cal G}$ and suppose that $S$ and $T$ are disjoint induced subgraphs of $G$ which together contain the vertex set of $G$. The following conditions are equivalent:

1. $S$ and $T$ are acyclic and $S\cap C$  is a path or a vertex  for every facial cycle~$C$ of $G$,

2. $ \{e^{*} \in E(G^{*}) : e \in  E[S, T] \}$ is  the edge set of a Hamilton cycle in $G^{*}$, 

3. $S$ and $T$ are trees,

4. $S$ is a  tree and $S\cap C$ is a path or a vertex for every facial cycle $C$ of $G$.
\end{theorem}
%$$\mbox{\bf Proof  of Theorem \ref{theorem1}}$$  
%\section{Proof of Theorem $1$}
\begin{PrfFact} 
Let $S$ and $T$ be disjoint induced subgraphs of $G$ which together contain the vertex set of $G$. We denote $E[S, T]:= E[V(S), V(T)]$. Notice that  $\{E(S), E(T), E[S,T]\}$ is a partition of $E(G)$.

Proof. $1 \Rightarrow 2$: Assume that $S$ and $T$ are acyclic and suppose that $S\cap C$  is a path  or a vertex for every facial cycle $C$ of $G$. Let $H$  denote a subgraph of $G^{*}$  such that $E(H) = (E[S, T])^{*}$.  Suppose that $H$  is not connected. Then  there exists an edge  cut of $G^{*}$ disjoint with $H$, because $G^{*}$ is connected. Thus, there exists a bond $B$ of $G^{*}$ disjoint with $H$. Hence, $B$ is a bond contained in $S^{*} \cup T^{*}$, because $\{S^{*}, T^{*}, (E[S, T])^{*}\}$) is a partition of $E(G^{*})$. Then, $B^{*}$ is a cycle contained in $S\cup T$. We obtain a contradiction, because $S \cup T$ is a forest. Hence, the graph $H$ is connected.

Let $v \in V(G^{*})$ and suppose that $C$ is a facial cycle of $v^{*}$ in $G$.  Then $C$ contains exactly two  edges of $E[S, T]$, because $S\cap C$  and $T\cap C$ are  paths or one-vertex subsets which together contain all vertices of $C$. Thus, $v$ has degree~$2$ in $H$. Hence, $H$ is a Hamilton cycle in $G^{*}$.

Proof. $2 \Rightarrow 3$ (see Bondy and  Murty \cite{flobar1} 10.2.11):  Assume that $(E[S, T])^{*}$ is  the edge set of a Hamilton cycle in $G^{*}$. Hence, $E[S, T]$ is a bond of $G$. Thus, both $S$ and $T$ are connected.

If $C$ is a cycle in $G$, then $C^{*}$ is a bond of $G^{*}$. Hence, $C^{*}$ contains at least two edges of the Hamilton cycle.  Thus, $C$ contains at least two edges belonging to $E[S, T]$ . Hence, neither $S$ nor $T$ contain $C$.  Therefore,  both $S$ and $T$ are acyclic. 

Proof. $3 \Rightarrow 4$: Assume that $S$ and $T$ are trees.  If $C$ is a facial cycle in $G$, then it has a common vertex with $S$ and with $T$, because $S$ and $T$ are acyclic. By the Jordan curve theorem,  $S\cap C$  is a path or a vertex, because $S$ and $T$ are disjoint induced trees which together contain all vertices of $C$.

Proof. $4 \Rightarrow 1$:  Assume that $S$ is a  tree such that $S\cap C$ is a path or a vertex for every facial cycle $C$ of $G$.  If $T$ contains a cycle (say $D$), then, by the Jordan curve theorem, $S$ lies either inside or outside $D$, since $S$ is connected. The other side of $D$ contains a face of $G$. Since the  boundary of this  face is disjoint with $S$, we obtain a contradiction.
\end{PrfFact}

\end{document}